\documentclass[11pt]{article} 
\usepackage{epic}
\usepackage{eepic}
\renewcommand{\theequation}{1.\arabic{equation}}
\begin{document}

\newcommand{\un}{\underline}
\author{ }
\newcommand{\no}{\noindent}
\newcommand{\be}{\begin{equation}}
\newcommand{\ee}{\end{equation}}
\newcommand{\ba}{\begin{array}}
\newcommand{\ea}{\end{array}}
\newcommand{\bea}{\begin{eqnarray*}}
\newcommand{\eea}{\end{eqnarray*}}
\newcommand{\bean}{\begin{eqnarray}}
\newcommand{\eean}{\end{eqnarray}}

 \newtheorem{theorem}{Theorem}
\newcommand{\bth}{\begin{theorem}}
\newtheorem{lemma}{Lemma}
   \newcommand{\blem}{\begin{lemma}}
   \newcommand{\elem}{\end{lemma}}
\newtheorem{corollary}{Corollary}
   \newcommand{\bcor}{\begin{corollary}}
   \newcommand{\ecor}{\end{corollary}}

\makeatletter
\@addtoreset{equation}{section}
\makeatother

{\centerline{{\Large {\bf A limiting form of the $q-$Dixon 
${}_4\varphi_3$ summation}}}}
{\centerline{\Large{\bf and related partition identities }}}
\centerline{by}
{\centerline{$Krishnaswami$ $Alladi$ $^*$ $and$ $Alexander$ $Berkovich$
{\footnote[1]{{Research supported in part by National Science Foundation 
Grant DMS-0088975}
\newline \hbox{}\qquad{2000 Mathematics Subject Classification: 05A17, 
05A19, 11P83, 11P81, 33D15, 33D20} 
\newline \hbox{}\qquad {\em Key words and phrases}: $q-$Dixon
Sum, $q-$Dougall Sum, weighted partition identities,
G\"ollnitz's (Big) \newline \hbox{}\qquad theorem, quartic reformulation, 
Jacobi's triple product, G\"ollnitz-Gordon identities, modular relations.}}\\}}

\abstract
By considering a limiting form of the $q-$Dixon ${}_4\varphi_3$ summation, 
we prove a weighted partition theorem involving odd parts differing by 
$\ge 4$. A two parameter refinement of this theorem is then deduced from 
a quartic reformulation of G\"ollnitz's (Big) theorem due to Alladi, and 
this leads to a two parameter extension of Jacobi's triple product identity 
for theta functions. Finally, refinements of certain modular identities of 
Alladi connected to the G\"ollnitz-Gordon series are shown to follow from 
a limiting form of the $q-$Dixon ${}_4\varphi_3$ summation.
\endabstract

\centerline {\bf{\S1: Introduction}}

The $q-$hypergeometric function $_{r+1}\varphi_r$ in $r+1$ numerator
parameters $a_1, a_2, \dots, a_{r+1},$ and $r$ denominator parameters
$b_1, b_2, \dots, b_r$, with base $q$ and variable $t$, is defined as
\be 
{}_{r+1}\varphi_r\left(\ba{cc}a_1, a_2, a_3, \dots, a_{r+1} \\ b_1, b_2,
\dots, b_r\ea ; q, t \right) =\sum^{\infty}_{k=0} \frac{(a_1;
q)_k\dots(a_{r+1}; q)_kt^k}{(b_1; q)_k(b_2;q)_k\dots(b_r ; q)_k(q; q)_k}.
\ee
When $r=3$, for certain special choices of parameters $a_i, b_j$, and
variable $t$, it is possible to evaluate the sum on the right in (1.1)
to be a product.  More precisely, the ${}_4\varphi_3$ $q-$Dixon summation
([4], (II. 13), p. 237) is 
\be 
{}_4\varphi_3\left(\ba{cc}a,-q\sqrt a, b, c \\ -\sqrt a,
\frac{aq}{b}, \frac{aq}{c}\ea; q, \frac{q\sqrt
a}{bc}\right)=\frac{(aq;q)_\infty (\frac{q\sqrt a}{b}; q)_\infty
(\frac{q\sqrt a}{c}; q)_\infty (\frac{aq}{bc};
q)_\infty}{(\frac{aq}{b}; q)_\infty (\frac{aq}{c};q)_\infty (q\sqrt a;
q)_\infty(\frac{q\sqrt a}{bc}; q)_\infty}.
\ee
Here and in what follows we have made use of the standard notation
\be 
(a;q)_n=(a)_n=\left\{\begin{array}{lll}
\prod^{n-1}_{j=0}(1-aq^j), & \mbox{ if } & n>0,\\ 
1, & \mbox{ if } & n=0,\\ 
\end{array} \right.
\ee
for any complex number $a$ and a non-negative integer $n$, and 
\be
(a)_\infty=(a;q)_\infty=\lim\limits_{n\to \infty} (a; q)_n
=\prod^\infty_{j=0} (1-aq^j),\textnormal{ for }|q|<1.
\ee
Sometimes, as in (1.3) and (1.4), when the base is $q$, we might
suppress it, but when the base is anything other than $q$, it will be
made explicit.

Our first goal is to prove Theorem 1 is \S2 which is a weighted identity
connecting partitions into odd parts differing by $\ge 4$ and
partitions into distinct parts $\not\equiv$ 2(mod 4).  We achieve this by
showing that the analytic representation of Theorem 1 is
\be
\sum_{k=0}^{\infty}\frac{z^k q^{4T_{k-1}}q^{3k}(z^2 q^2; q^2)_k
(1+zq^{2k+1})}{(q^2; q^2)_k}=(-zq;q^2)_\infty(z^2q^4; q^4)_\infty
\ee
and establish (1.5) by utilizing a limiting form of the $q-$Dixon
summation formula (1.2). In (1.5), $T_k=k(k+1)/2$ is the $k-$th
triangular number.

It is possible to obtain a two parameter refinement of Theorem 1 by
splitting the odd integers into residue classes 1 and 3 (mod 4) and
keeping track of the number of parts in each of these residue
classes.  This result, which is stated as Theorem 2 in \S 3, is a
special case of a weighted reformulation of G\"ollnitz's (Big) theorem
due to Alladi (Theorem 6 of [2]). In \S 3 we also state an analytic
identity (see (3.3)) in two free parameters $a$ and $b$ that is
equivalent to Theorem 2, and note that (1.5) follows from this as the 
special case $a=b=z$.  Identity (3.3) can be viewed as a two parameter 
generalization of Jacobi's celebrated triple product identity for 
theta functions (see (3.5) in \S 3). 

Identity (3.3) is itself a special case of {\it{key identity}} in three 
free parameters $a, b,$ and $c$, due to Alladi and Andrews ([3], eqn. 3.14), 
for G\"ollnitz's (Big) Theorem.  The proof of this key identity of Alladi 
and Andrews in [3] utilizes Jackson's $q-$analog of Dougall's summation 
for ${}_6\varphi_5$. Note that the left hand side of (3.3) is a double 
summation.  On the other hand, the left side of (1.5) is just a single 
summation, and its proof requires only a limiting form of the $q-$Dixon 
summation for ${}_4\varphi_3$. Owing to the choice $a=b=z$, the double 
sum in (3.3) reduces to a single summation in (4.6) resembling (1.5), 
and this process is described in \S4. Finally, certain modular identities 
for G\"ollnitz-Gordon functions due to Alladi [2] are refined in \S5 
using a limiting case of the $q-$Dixon summation (1.2).

We conclude this section by mentioning some notation pertaining to
partitions.  For a partition $\pi$ we let

\centerline{$\sigma(\pi)=$ the sum of all parts of $\pi$,}
\centerline{$\nu(\pi)=$ the number of parts of $\pi$,}
\centerline{$\nu(\pi; r, m)=$ the number of parts of $\pi$ which are
$\equiv r$(mod $m$),}
\centerline {$\nu_d(m)=$ the number of different parts of $\pi$,}
\centerline{$\nu_{d, \ell}(m)=$ the number of different parts of $\pi$ 
which are $\ge\ell$, and} 
\centerline{$\lambda (\pi)=$ the least part of $\pi$.}

\newpage

\renewcommand{\theequation}{2.\arabic{equation}}
\setcounter{equation}{0}

\centerline{\bf{\S2: Combinatorial interpretation and proof of (1.5)}}
\bigskip
Let ${\cal O}_4$ denote the set of partitions into odd parts differing
by $\ge 4$.  Given $\tilde\pi \in {\cal O}_4$, a chain $\chi$ in
$\tilde\pi$ is defined to be a maximal string of consecutive parts
differing by exactly 4.  Let $N_\lambda(\tilde\pi)$ denote the number
of chains in $\tilde\pi$ with least part $\ge \lambda$.

Next let ${\cal D}_{2, 4}$ denote the set of partitions into distinct parts
$\not\equiv 2(\textnormal{mod}$ $4)$. We then have

{\bf{Theorem 1}}: {\it{For all integers n}}$\ge 0$ 
$$\sum_{\ba{cc} \tilde\pi \in {{\cal O}_4} \\ \sigma(\tilde\pi)=n\ea} 
z^{\nu(\tilde\pi)}(1-z^2)^{N_5(\tilde\pi)}
=\sum_{\ba{cc} \pi \in {\cal D}_{2,4} \\ \sigma(\pi)=n\ea} 
z^{\nu(\pi;1,2)}(-z^2)^{\nu(\pi;0,2)}
$$

So, for partitions $\tilde\pi\in {\cal O}_4$, we attach the weight $z$
to each part, and the weight $(1-z^2)$ to each chain having least part
$\ge 5$.  The weight of $\tilde\pi$ is then defined multiplicatively.
Similarly, for partitions $\pi\in {\cal D}_{2, 4}$, each odd part is assigned
weight $z$, and each even part is assigned the weight $-z^2$, where all these 
even parts are actually multiples of 4. For example, when $n=10$, the 
partitions in ${\cal O}_4$ are $9+1$ and $7+3$, with weights
$z^2(1-z^2)$ and $z^2$ respectively. These weights add up to yield
$2z^2-z^4$.  The partitions of 10 in ${\cal D}_{2, 4}$ are $9+1, 7+3,$ and
$5+4+1$, with weights $z^2, z^2$, and $z^2(-z^2)$ respectively.  These
weights also add up to $2z^2-z^4$, verifying Theorem 1 for $n=10$.

We will now show that Theorem 1 is the combinational interpretation of
(1.5).

It is clear that the product
\be
\prod^{\infty}_{m=1}(1+zq^{2m-1})(1-z^2 q^{4m})
\ee
on the right in (1.5) is the generating function of partitions $\pi\in
{\cal D}_{2, 4}$, with weights as in Theorem 1.  So we need to show that the
series on the left in (1.5) is the generating function of partitions
$\tilde\pi\in{\cal O}_4$ with weights as specified in Theorem 1.  For
this we consider two cases.

\un{Case 1}:  $\lambda(\tilde\pi)\ne 1$.  

If $\tilde\pi$ is non-empty, then $\lambda(\tilde\pi)\ge 3$.

Since the parts of $\tilde\pi$ differ by $\ge 4$, we may subtract $0$
from the smallest part, 4 from the second smallest part, 8 from the
third smallest, ..., $4k-4$ from the largest part of $\tilde\pi$,
assuming $\nu(\tilde\pi)=k$.  We call this procedure the {\it{Euler
subtraction}}. After the Euler subtraction is performed on $\tilde\pi$,
we are left with a partition $\pi'$ into $k$ odd parts such that the
number of different parts of $\pi'$ is precisely the number of chains
in $\tilde\pi$.  If we denote by $G_{3, k}(q, z)$ the generating
function of partitions $\tilde\pi\in{\cal O}_4$ with $\lambda(\tilde\pi)\ne
1, \nu(\tilde\pi)=k$, and counted with weight $z^{\nu(\tilde\pi)}
(1-z^2)^{N_5(\tilde\pi)}$, then the Euler subtraction process yields
\be 
G_{3, k}(q, z)=z^k q^{4T_{k-1}} g_{3, k} (q, z),
\ee
where $g_{3, k}(q, z)$ is the generating function of partitions $\pi'$
into $k$ odd parts each $\ge 3$ and counted with weight
$(1-z^2)^{\nu_{d, 5}(\pi ')}$.  

At this stage we make the observation that if a set of positive integers J 
is given, then   
\be
\prod_{j\in J}\left (\frac{1-twq^j}{1-tq^j}\right)
=\prod_{j\in J}(1+(1-w)\left\{tq^j+t^2 q^{2j}+t^3 q^{3j}+\dots\right\})
\ee
is the generating function of partitions $\pi^*$ into parts belonging to $J$ 
and counted with weight $t^{\nu(\pi^*)}(1-w)^{v_d(\pi^\ast)}$. So from the 
principles underlying (2.3) it follows that 
\be
\sum^{\infty}_{k=0}g_{3, k}(q,z)t^k
=\frac{1}{(1-tq^3)}\prod^{\infty}_{j=0}\left(\frac{1-tz^2q^{2j+5}}
{1-tq^{2j+5}}\right)=\frac{(tz^2q^5; q^2)_\infty}{(tq^3; q^2)_\infty}.
\ee
Using Cauchy's identity
\be
\frac{(at)_\infty}{(t)_\infty}=\sum^{\infty}_{k=0}\frac{(a)_kt^k}{(q)_k},
\ee
we can expand the product on the right in (2.4) as
\be
\frac{(tz^2q^5; q^2)_\infty}{(tq^3;
q^2)_\infty}=\sum^{\infty}_{k=0}\frac{t^kq^{3k}(z^2q^2; q^2)_k}{(q^2;q^2)_k}.
\ee
So by comparing the coefficients of $t^k$ in (2.4) and (2.6) we get
\be 
g_{3, k}(q;z)=\frac{q^{3k}(z^2q^2;q^2)_k}{(q^2;q^2)_k}.
\ee
Thus (2.7) and (2.2) yield
\be 
G_{3, k}(q;z)=\frac{z^k q^{3k} q^{4T_{k-1}}(z^2 q^2; q^2)_k}{(q^2;
q^2)_k}, \textnormal{ for } k\ge 0.
\ee

\un{Case 2}:  $\lambda(\tilde\pi)=1$.

Here for $k>0$ we denote by $G^\ast_{1, k}(q, z)$ the generating function of
partitions $\tilde\pi\in{\cal O}_4$ having $\lambda(\tilde\pi)=1$, 
$\nu(\tilde\pi)=k$, and  counted with weight
$z^{\nu(\tilde\pi)}(1-z^2)^{N_5(\tilde\pi)}$. The Euler subtraction
process yields
\be 
G^\ast_{1, k}(q, z)=z^k q^{4T_{k-1}}g^\ast_{1, k}(q, z),
\ee
where $g^\ast_{1, k}(q, z)$ is the generating function of partitions
$\pi'$ into $k$ odd parts counted with weight $(1-z^2)^{\nu_{d,3}(\pi ')}$.  
The principles underlying (2.3) show that 
\be
\sum^{\infty}_{j=1}g^*_1(q, z) t^k=\frac{tq}{1-tq}\prod^{\infty}_{j=1}
\left(\frac{1-tz^2q^{2j+1}}{1-tq^{2j+1}}\right)
=tq\frac{(tz^2q^3; q^2)_\infty}{(tq; q^2)_\infty}
=\sum^{\infty}_{k=0} t^{k+1}q^{k+1}\frac{(z^2 q^2; q^2)_k}{(q^2;q^2)_k}
\ee
by Cauchy's identity (2.5).  Thus by comparing the coefficients of
$t^k$ at the extreme ends of (2.10), we get
\be 
g^\ast_{1, k}(q, z)= q^k\frac{(z^2 q^2; q^2)_{k-1}}{(q^2;q^2)_{k-1}}.
\ee
This when combined with (2.9) yields
\be 
G^*_{1, k}(q, z)=\frac{z^kq^kq^{4T_{k-1}}(z^2q^2:q^2)_{k-1}}{(q^2;q^2)_{k-1}}.
\ee

Finally, it is clear that
\be
\sum^{\infty}_{k=0}G_{3, k}(q, z)+\sum^{\infty}_{k=1}G^\ast_{1,
k}(q, z)
\ee
is the generating function of partitions $\tilde\pi\in{\cal O}_4$
counted with weight as specified in Theorem 1. From (2.8) and (2.12),
the sum in (2.13) can be seen to be
\bean
\sum^{\infty}_{k=0} z^k q^{3k}q^{4T_{k-1}}
\frac{(z^2 q^2;q^2)_k}{(q^2;q^2)_k}+\sum^{\infty}_{k=1}z^kq^kq^{4T_{k-1}}
\frac{(z^2q^2;q^2)_{k-1}}{(q^2; q^2)_{k-1}}\nonumber \\
=\sum^{\infty}_{k=0} z^k q^{3k} q^{4T_{k-1}}\frac{(z^2 q^2;
q^2)_k}{(q^2;q^2)_k}+\sum^{\infty}_{k=0}
z^{k+1}q^{k+1}q^{4T_k}\frac{(z^2 q^2; q^2)_k}{(q^2; q^2)_k}\nonumber \\ 
=\sum^{\infty}_{k=0}\frac{z^kq^{3k}q^{4T_{k-1}}(z^2q^2;q^2)_k(1+zq^{2k+1})}{(q^2; q^2)_k}
\eean
which is the series on the left in (1.5).

From (2.14) and (2.1) it follows that Theorem 1 is the combinatorial
interpretation of (1.5).  Thus to prove Theorem 1 it suffices to
establish (1.5) and this what we do next.

From the definition of the $q-$hypergeometric function in (1.1) we see that 
\bean
(1+zq)
{}_4\varphi_3\left(\ba{cc}z^2q^2,-zq^3, \rho, \rho \\ 
-zq,\frac{z^2q^4}{\rho}, \frac{z^2q^4}{\rho}\ea; q^2, \frac{zq^3}{\rho^2}
\right) \nonumber \\ 
=(1+zq)\sum^{\infty}_{k=0}\frac{(z^2 q^2; q^2)_k(-zq^3; q^2)_k
(\rho;q^2)_k(\rho; q^2)_k}
{(q^2; q^2)_k(-zq; q^2)_k(\frac{z^2q^4}{\rho}; q^2)_k
(\frac{z^2q^4}{\rho}; q^2)_k}\frac{z^kq^{3k}}{\rho^{2k}}
\nonumber \\ \nonumber \\ 
=\sum^{\infty}_{k=0}\frac{z^kq^{3k}(z^2q^2;q^2)_k(1+zq^{2k+1})}{(q^2;q^2)_k}
\frac{(\rho;q^2)_k(\rho;q^2)_k}{(\frac{z^2q^4}{\rho};q^2)_k
(\frac{z^2q^4}{\rho};q^2)_k}\frac{1}{\rho^{2k}}.
\eean

Next observe that 
\be
\lim_{\rho\to\infty}\frac{(\rho;q^2)_k}{\rho^k}
=\lim_{\rho\to\infty}\frac{(1-\rho)(1-\rho q^2)\dots(1-\rho q^{2k-2})}{\rho^k}
=(-1)^kq^{2T_{k-1}}
\ee
and
\be
\lim_{\rho\to\infty}\left(\frac{z^2q^4}{\rho}; q^2\right)_k =1.
\ee
Thus (2.15), (2.16), and (2.17) imply that
\be
(1+zq)\lim_{\rho\to\infty} 
{}_4\varphi_3\left(\ba{cc}z^2q^2, -zq^3, \rho, \rho \\ 
-zq, \frac{z^2q^4}{\rho}, \frac{z^2q^4}{\rho}\ea;q^2, \frac{zq^3}{\rho^2}
\right)
=\sum^{\infty}_{k=0}\frac{z^kq^{3k}q^{4T_{k-1}}(z^2q^2; q^2)_k(1+zq^{2k+1})}
{(q^2; q^2)_k}
\ee
which is the series on the left in (1.5).

At this stage we observe that the hypergeometric function ${}_4\varphi_3$
on the left in (2.18) is precisely the one in the $q-$Dixon summation
(1.2) with the replacements
\be 
q\mapsto q^2, a\mapsto z^2q^2, b\mapsto \rho, c\mapsto \rho.
\ee

Thus with substitutions (2.19) in (1.2) we deduce that 
\bean
(1+zq)\lim_{\rho\to \infty} 
{}_4\varphi_3\left(\ba{cc}z^2q^2, -zq^3, \rho, \rho \\ 
-zq, \frac{z^2q^4}{\rho}, \frac{z^2q^4}{\rho}\ea ; q^2, \frac{zq^3}{\rho^2}
\right)\nonumber \\ \nonumber \\
=(1+zq)\lim_{\rho\to \infty}\frac{(z^2q^4; q^2)_\infty 
(\frac{zq^3}{\rho}; q^2)_\infty
(\frac{zq^3}{\rho}; q^2)_\infty (\frac{z^2q^4}{\rho^2};q^2)_\infty}
{(\frac{z^2q^4}{\rho}; q^2)_\infty (\frac{z^2q^4}{\rho};q^2)_\infty
(zq^3;q^2)_\infty (\frac{zq^3}{\rho^2}; q^2)_\infty}
\nonumber \\ \nonumber \\
=(1+zq)\frac{(z^2q^4; q^2)_\infty}{(zq^3; q^2)_\infty}
=\frac{(z^2q^2;q^2)_\infty}{(zq;q^2)_\infty}
\nonumber \\ 
=\frac{(z^2q^2; q^4)_\infty(z^2q^4; q^4)_\infty}{(zq;q^2)_\infty}
=(-zq;q^2)_\infty(z^2q^4; q^4)_\infty.
\eean
Thus (1.5) follows from (2.18) and (2.20) and this completes the proof
of Theorem 1.
\bigskip
\renewcommand{\theequation}{3.\arabic{equation}}
\setcounter{equation}{0}

\centerline{\bf{\S3: Two parameter refinement}}
\bigskip
When a partition $\tilde\pi\in{\cal O}_4$ is decomposed into chains, the
parts in a given chain all belong to the same residue class mod 4.
This suggests that there ought to be a two parameter refinement of
Theorem 1 in which we can keep track of parts in residue classes 1 and
3 (mod 4) separately.  Theorem 2 stated below is such a
refinement. Actually Theorem 2 is a special case of a refinement and
reformulation of a deep theorem of G\"ollnitz [5] in three parameters
$a,b$, and $c$ due to Alladi ([2], Theorem 6) by setting one of the
parameters $c=-ab$.

{\bf{Theorem 2}:} {\it{For all integers $n\ge 0$ and complex numbers 
$a$ and $b$ we have}}
$$
\sum_{\ba{cc}_{\tilde\pi\in{\cal O}_4} \nonumber \\ 
\sigma(\tilde\pi)=n\ea}
a^{\nu(\tilde\pi; 1, 4)}b^{\nu(\tilde\pi;3,4)}(1-ab)^{N_5(\tilde\pi)} 
=\sum_{\ba{cc}\pi\in {\cal D}_{2, 4}\\ \sigma(\pi)=n\ea} 
a^{\nu(\pi; 1, 4)}b^{\nu(\pi;3,4)} (-ab)^{\nu(\pi,0,4)}
$$

Since Theorem 2 is a two parameter refinement of Theorem 1 which has 
the analytic representation (1.5), it is natural to ask for an analytic 
identity in two free parameters that reduces to (1.5). We will now obtain 
such a two parameter identity, namely, (3.3) below. Instead of deriving 
(3.3) combinatorially from Theorem 2 by following the method in \S2, we 
will now illustrate a different approach which involves a certain cubic 
reformulation of G\"ollnitz's (Big) theorem due to Alladi [1], and its 
{\it{key identity}} in three free parameters $a$, $b$, and $c$ due to 
Alladi and Andrews [3]. More precisely, we will show that (3.3) is the 
analytic representation of Theorem 2 after a discussion of the 
following special case $c=-ab$ of the key identity (3.14) of [3]:  
$$ 
\sum_{i,j\ge 0}\frac{a^i q^{(3i^2-i)/2}(ab; q^3)_i}{(q^3;q^3)_i}\cdot 
\frac{b^j q^{3ij} q^{(3j^2+j)/2}(ab;q^3)_j}{(q^3; q^3)_j}\cdot 
\left(\frac{1-abq^{3(i+j)}}{1-ab}\right)
$$
\be
=(-aq;q^3)_\infty (-bq^2; q^3)_\infty(abq^3; q^3)_\infty.
\ee

The cubic reformulation of G\"ollnitz's theorem in [1] in three free 
parameters $a$, $b$, and $c$, was in the form of an identity connecting 
partitions into distinct parts with a weighted count of partitions into 
parts differing by $\ge 3$. When we set $c=-ab$ in Theorem 2 of [1], 
it turns out that within the set of partitions into parts differing by 
$\ge 3$, we need only consider those partitions not having any multiples of 3 
as parts; this is because (see [1], Theorem 2) the choice $c=-ab$ 
makes the weights equal to 0 if the partition has a multiple of 3 in it. 
Thus from the analysis in [3] and the specialization $c=-ab$ in Theorem 2 
of [1], it follows that the combinatorial interpretation of (3.1) is

{\bf{Theorem 3}:} {\it{Let ${\cal D}$ denote the set of partitions into 
distinct parts. Let ${\cal D}_3$ denote the set of partitions into parts 
differing by $\ge 3,$ and containing no multiples of 3. Given 
$\tilde\pi\in {\cal D}_3$, decompose it into chains, where a chain here 
is a maximal string of parts differing by exactly 3.  Then we have}}
$$
\sum_{\ba{cc}_{\tilde\pi\in {\cal D}_3} \\ \sigma(\tilde\pi)=n\ea} 
a^{\nu(\tilde\pi;1,3)}b^{\nu(\tilde\pi;2,3)}(1-ab)^{N_3(\tilde\pi)}
=\sum_{\ba{cc} \pi\in {\cal D} \\ \sigma(\pi)=n\ea} 
a^{\nu(\pi;1,3)}b^{(\pi;2,3)} (-ab)^{\nu(\pi; 0,3)}.
$$

As an example, when $n=9$, the partitions in ${\cal D}_3$ are
8+1 and 7+2 both having weights $ab(1-ab)$.  So these weights will add
up to $2ab(1-ab)$.  The partitions of 9 in ${\cal D}$ with their
corresponding weights are listed below:

\un{Partitions}: $9, \quad 8+1, \quad 7+2, \quad 6+3, \quad 6+2+1,
\quad 5+4, \quad 5+3+1, \quad 4+3+2$

\un{Weights}:  $-ab, \quad ab, \quad\quad\quad ab, \quad (-ab)^2, 
\quad ab(-ab), \quad ab, \quad\quad ab(-ab), \quad ab(-ab)$

The above weights when added also yield $2ab(1-ab)$ thereby verifying
Theorem 3 for $n=9$.

Theorems 3 and 2 are really the same because the role of the modulus 3 in 
Theorem 3 is replaced by the modulus 4 in Theorem 2.  More precisely, we may 
view ${\cal D}_3$ and ${\cal O}_4$ as sets of partitions into parts differing 
by $\ge k$ and containing only parts in the residue classes 
$\pm 1(mod\quad k)$, for $k=3,4$. Similarly, we may think of $\cal D$ and 
${\cal D}_{2, 4}$ as sets of partitions into parts 
$\equiv 0,\pm 1(mod\quad k)$, 
for $k=3,4$. Note that the functions $\nu(\pi;r,3)$ and $\nu(\tilde\pi;r,3)$ 
in Theorem 3 are replaced $\nu(\pi;r,4)$ and $\nu(\tilde\pi;r,4)$ in 
Theorem 2. Pursuing this line of correspondence, $N_3(\tilde\pi)$ in 
Theorem 3 should be replaced by $N_4(\tilde\pi)$ in Theorem 2, but this 
is the same as having $N_5(\tilde\pi)$ in Theorem 2 because 
$\tilde\pi\in{\cal O}_4$.

Having observed the correspondence between Theorems 2 and 3, we deduce
that the analytic representation of Theorem 2 ( in the form of a two
parameter $q-$hypergeometric identity) is obtained from (3.1) by the
substitutions 
\be 
q\mapsto q^{4/3}, a \mapsto aq^{-1/3}, b\mapsto bq^{1/3}, ab\mapsto ab,
\ee
which yields
$$ 
\sum_{i,j\ge0}\frac{a^i q^{2i^2-i}(ab;q^4)_i}{(q^4; q^4)_i}\cdot
\frac{b^jq^{4ij}q^{2j^2+j}(ab; q^4)_j}{(q^4; q^4)_j}\cdot
\left(\frac{1-abq^{4(i+j)}}{1-ab}\right)
$$
\be
=(-aq; q^4)_\infty (-bq^3; q^4)_\infty (ab;q^4; q^4)_\infty.
\ee

Identities (3.1) and (3.3) are interesting for another reason.  They
can be considered as two parameter extensions of Jacobi's celebrated
triple product identity for theta functions.  More precisely, if we
put $ab=1$ in (3.1) and (3.3), then on the left hand side of each of
these identities, only the terms having either $i=0$ or $j=0$ survive,
and so the identities reduce to 
\be
\sum^{\infty}_{i=-\infty}a^i q^{(3i^2 -i)/2}
=(-aq; q^3)_\infty(-a^{-1}q^2;q^3)_\infty(q^3; q^3)_\infty
\ee
and
\be
\sum^{\infty}_{i=-\infty}a^i q^{2i^2-i}
=(-aq; q^4)_\infty(-a^{-1}q^3;q^4)_\infty (q^4;q^4)_\infty
\ee
which are equivalent to Jacobi's triple product identity.

Identities (3.4) and (3.5) can be deduced combinatorially from
Theorems 3 and 2 respectively, by setting $ab=1$.  This is because
$ab=1$ forces $N_3(\tilde\pi)=0$ (resp. $N_5(\tilde\pi)=0)$ in Theorem
3 (resp. Theorem 2) and this brings about a drastic reduction in the
type of partitions to be enumerated in ${\cal D}_3$ (resp. ${\cal O}_4$). We 
refer the reader to Alladi [1], [2] for these combinatorial arguments.
\bigskip

\renewcommand{\theequation}{4.\arabic{equation}}
\setcounter{equation}{0}

\centerline{\bf{\S4: Reduction to a single summation}}
\bigskip
The left hand side of identity (3.3) is a double summation.  It turns out 
that if we set
$$
a=b=z,
$$
then the left side of (3.3) reduces to a single infinite sum.  It is quite 
instructive to see how this happens, and so we describe it now.

First observe that
\be 
2i^2-i+2j^2+j+4ij=2(i+j)^2-(i+j)+2j.
\ee
Thus if we set $a=b=z$ and reassemble the terms in (3.3) with $k=i+j$,
then (4.1) shows that (3.3) becomes
\be
\sum_{k=0}^\infty z^kq^{2k^2-k}\frac{(1-z^2q^{4k})}{(1-z^2)}
\left\{\sum_{i+j=k}\frac{q^{2j}(z^2;q^4)_i(z^2; q^4)_j}
{(q^4; q^4)_i(q^4; q^4)_j}\right\}
=(-zq;q^2)_\infty(z^2 q^4; q^4)_\infty.
\ee
At this point we take the product in Cauchy's identity (2.5) and
decompose it as 
\be 
\frac{(at)_\infty}{(t)_\infty}
=\frac{(at; q^2)_\infty}{(t;q^2)_\infty}\cdot
\frac{(atq; q^2)_\infty}{(tq; q^2)_\infty}.
\ee
If we now substitute the expansion in (2.5) for each of the products
in (4.3), we get
\be 
\sum^{\infty}_{k=0}\frac{(a)_kt^k}{(q)_k}
=\left(\sum^\infty_{i=0}\frac{(a;q^2)_it^i}{(q^2;q^2)_i}\right)
\left(\sum^\infty_{j=0}\frac{(a;q^2)_jt^jq^j}{(q^2;q^2)_j}\right).
\ee
By comparing the coefficients of $t^k$ on both sides of (4.4) we obtain
\be
\frac{(a)_k}{(q)_k}=\sum_{i+j=k}\frac{(a;q^2)_i(a;q^2)_j q^j}
{(q^2;q^2)_i(q^2; q^2)_j}.
\ee
If we replace $q\mapsto q^2$ and $a\mapsto z^2$, we see that the sum in
(4.5) becomes the expression within the parenthesis (namely the inner
sum) on the left in (4.2).  Thus with these replacements (4.5) implies
that (4.2) can be written as the single summation identity
\be 
1+\sum^{\infty}_{k=1}\frac{z^k q^{2k^2-k}(z^2q^2; q^2)_{k-1}(1-z^2q^{4k})}
{(q^2; q^2)_k}=(-zq;q^2)_\infty(z^2 q^4; q^4)_\infty.
\ee

Identity (4.6) is an analytic representation of Theorem 1.  Note
however that the series in (4.6) is different from the series in
(1.5).  The explanation of this is as follows.

If we add $G_{3, k}(q, z)$ and $G^*_{1,k}(q, z)$ for each $k\ge 1$,
we get from (2.8) and (2.12)
$$
G_{3, k}(q, z)+G^*_{1, k}(q, z)
=\frac{z^kq^k q^{4T_{k-1}}(z^2 q^2;q^2)_{k-1}}{(q^2; q^2)_k}
\left\{1-q^{2k}+q^{2k}(1-z^2 q^{2k})\right\}
$$
$$
=\frac{z^kq^{2k^2-k}(z^2q^2;q^2)_{k-1}(1-z^2q^{4k})}{(q^2;q^2)_k}
$$
which is the $k-$th summand in (4.6). The starting term 1 in (4.6) is
to be interpreted as $G_{3,0}(q, z)$. On the other hand in (2.14) we
are considering 
$$
G_{3, k}(q,z)+G^\ast_{1, k+1}(q,z),\textnormal{ for }k\ge 0,
$$
and this leads to the series in (1.5) which is different from (4.6).

The reason we preferred (1.5) to (4.6) is because (1.5) could be
proved using only a limiting form of the $q-$Dixon summation of
${}_4\varphi_3$, whereas (4.6) would have required a limiting form of Jackson's
${}_6\varphi_5$ summation ([4], (II.21), p. 238).
\bigskip
\renewcommand{\theequation}{5.\arabic{equation}}
\setcounter{equation}{0}

\centerline{\bf{\S5: Modular identities for the refined G\"ollnitz-Gordon functions.}}
\bigskip
The well known G\"ollnitz-Gordon identities are
\be
G(q)=\sum^{\infty}_{n=0}\frac{q^{n^2}(-q; q^2)_n}{(q^2;q^2)_n}
=\frac{1}{(q; q^8)_\infty(q^4; q^8)_\infty (q^7; q^8)_\infty}
\ee
and
\be
H(q)=\sum^{\infty}_{n=0}\frac{q^{n^2+2n}(-q; q^2)_n}{(q^2;q^2)_n}
=\frac{1}{(q^3; q^8)_\infty (q^4; q^8)_\infty(q^5; q^8)_\infty}.
\ee
Identities (5.1) and (5.2) are actually (36) and (34) on Slater's list
[7], but it was G\"ollnitz [5] and Gordon [6] who realized their
partition significance and their relationship with a continued
fraction.  More precisely, the G\"ollnitz-Gordon partition theorem is:
{\it{For i=1,2, the number of partitions of an integer n into parts
differing by $\ge 2,$ with strict inequality if a part is even, and
least part $\ge 2i-1$, equals the number of partitions of n into
parts $\equiv 4, \pm (2i-1)$ (mod 8).}}

In view of the form of the series-product identities (5.1) and (5.2),
their partition interpretation given above, and their relationship with 
a certain continued fraction, the G\"ollnitz-Gordon identities are
considered as the perfect analogues for the modulus 8 for what the
celebrated Rogers-Ramanujan identities are for the modulus 5.

In [2] Alladi established four reformulations of G\"ollnitz's (Big)
partition theorem using four quartic transformations.  One of the
reformulations yielded Theorem 2 of \S3. Using one of the other
reformulations, Alladi [2] deduced the modular identity
\be 
G(-q^2)+qH(-q^2)=(-q;q^4)_\infty(q^2; q^4)_\infty(-q^3; q^4)_\infty
\ee
combinatorially.  Alladi [2] then defined the {\it{twisted}} 
G\"ollnitz-Gordon functions
\be 
G_t(q)=\sum^{\infty}_{n=0}\frac{q^{n^2}(q;q^2)_n}{(q^2; q^2)_n},
\ee
and 
\be
H_t(q)=\sum^{\infty}_{n=0}\frac{q^{n^2+2n}(q;q^2)_n}{(q^2;q^2)_n},
\ee
and deduced the modular identity
\be
G_t(q^2)+q H_t(q^2)=(-q; q^4)_\infty (-q^2; q^4)_\infty (q^3;q^4)_\infty
\ee
combinatorially from the same reformulation of G\"ollnitz's (Big)
theorem.

The twisted G\"ollnitz-Gordon functions do not have the product
representations of the type $G(q)$ and $H(g)$ possess.  But the modular
identity (5.6) implies that
\be
G_t(q^2)=\frac{(-q;q^4)_\infty(-q^2;q^4)_\infty (q^3; q^4)_\infty
+(q;q^4)_\infty(-q^2;q^4)_\infty(-q^3; q^4)_\infty}{2}
\ee
and 
\be
H_t(q^2)=\frac{(-q;q^4)_\infty(-q^2;q^4)_\infty (q^3;q^4)_\infty 
-(q;q^4)_\infty(-q^2;q^4)_\infty(q^3; q^4)_\infty}{2q}.
\ee
In the absence of product reprentations, (5.7) and (5.8) show that
$G_t(q^2)$ and $G_t(q^2)$ are arithmetic means of interesting
products.  From (5.3) it follows that $G(-q^2)$ and $H(-q^2)$ have
representations similar to (5.7) and (5.8).

By utilizing a limiting form of the $q-$Dixon summation (1.2) for
${}_4\varphi_3$, we will now establish a more general modular identity
(see (5.13) below) that contains both (5.3) and (5.6).  To this end
let $c=\delta\sqrt{aq}$ with $\delta=\pm 1$ in (1.2), and multiply
both sides by $1+\sqrt a$.  This way we find
\be
\sum^{\infty}_{k=0}\frac{\delta^kb^{-k}q^{k/2}(a)_k(b)_k(1+\sqrt a q^k)}{(q)_k 
\big(\frac{aq}{b}\big)_k}
=\frac{(a)_\infty(\frac{q\sqrt a}{b})_\infty (\delta\sqrt q)_\infty
(\frac{\delta\sqrt{aq}}{b})_\infty}
{(\frac{aq}{b})_\infty(\delta\sqrt{aq})_\infty(\sqrt a)_\infty
(\frac{\delta\sqrt q}{b})_\infty}.
\ee
Analogous to (2.16) we now have
\be
\lim_{b\rightarrow\infty}\frac{(b)_k}{b^k}=(-1)^k q^{T_{k-1}}.
\ee
Thus (5.9) and (5.10) imply that by going to the limit $b\rightarrow\infty$
we get
\be
\sum^{\infty}_{k=0}\frac{(-\delta)^k q^{{k^2}/2}(a)_k(1+\sqrt a q^k)}{(q)_k}
=\frac{(a)_\infty(\delta\sqrt q)_\infty}{(\sqrt a)_\infty
(\delta\sqrt{aq})_\infty}, \textnormal{ for }\delta=\pm 1.
\ee
In (5.11) if we replace $q\mapsto q^4$ and $a\mapsto a^2 q^2$, we obtain 
$$
\sum^{\infty}_{k=0}(-\delta)^k q^{2k^2}\frac{(a^2 q^2; q^4)_k}{(q^4;
q^4)_k}(1+aq^{4k+1})=\frac{(a^2 q^2; q^4)_\infty(\delta q^2;
q^4)_\infty}{(aq; q^4)_\infty(\delta aq^3; q^4)_\infty}
$$
\be
=\frac{(a^2q^2; q^8)_\infty(a^2 q^6; q^8)_\infty(\delta q^2;
q^4)_\infty}{(aq; q^4)_\infty (\delta aq^3; q^4)_\infty}=(-aq;
q^4)_\infty(\delta q^2; q^4)_\infty(-\delta aq^3; q^4)_\infty.
\ee
Note that the special case $\delta=1, a=1$ in (5.12) yields (5.3),
whereas $\delta=-1, a=1$ in (5.12) is (5.6).  We may write (5.12) in
the form of the modular identity
\be 
G_{a^2, \delta}(-q^2)+aq H_{a^2,\delta}(-q^2)
=(-aq;q^4)_\infty(\delta q^2; q^4)_\infty(-\delta aq^3;q^4)_\infty, 
\textnormal{ for }\delta=\pm 1
\ee
for the refined G\"ollnitz-Gordon functions which we define as
\be
G_{a,\delta}(q)=\sum^\infty_{k=0}\frac{\delta^k q^{k^2}(-aq; q^2)_k}
{(q^2;q^2)_k}
\ee
and
\be 
H_{a,\delta}(q)=\sum^{\infty}_{k=0}\frac{\delta^kq^{k^2+2k}(-aq;q^2)_k}
{(q^2; q^2)_k}.
\ee
Here $G_{a, \delta}(q)$ is the generating function of partitions into
parts differing by $\ge 2,$ with strict inequality if a part is even,
where each odd part is assigned weight $\delta$, and each even part
given weight $\delta a$.  The weight of the partition under
consideration is the product of the weights of its parts.  The
function $H_{a, \delta}(q)$ has a similar interpretation, except for
the added restriction that the least part in $\ge 3$ for the
partitions enumerated.

We note that the combinatorial arguments in [2] which yielded (5.3)
and (5.6) could be used to derive the more general modular relation (5.13).

\begin{center}
{\bf References}
\end{center}
\parskip=0pt
\baselineskip=12pt
\begin{enumerate}
\item[1.] \un{K. Alladi}, A combinatorial correspondence related to
G\"ollnitz's (Big) partition theorem and applications,
{\em Trans. Amer. Math Soc.} {\bf{349}} (1997), 2721-35.

\item[2.] \un{K. Alladi}, On a partition theorem of G\"ollnitz and
quartic transformations (with an appendix by B. Gordon), 
{\em J. Num. Th.} {\bf{69}} (1998), 153-180.

\item[3.] \un{K.  Alladi} and \un{G. E. Andrews}, A new key identity for
G\"ollnitz's (Big) partition theorem, {\em Contemp. Math.} {\bf{210}} (1998),
229-241.

\item[4.] \un{G. Gasper} and \un{M. Rahman}, Basic hyper-geometric series,
{\em Encyclopedia of Mathematics and its Applications}, Vol.{\bf{10}}, 
Cambridge (1990).

\item[5.] \un{ H. G\"ollnitz}, Partitionen mit Differenzenbedingungen,
{\em J. Reine Angew Math.} {\bf{225}} (1967), 154-190.

\item[6.] \un{B. Gordon}, Some continued fractions of the
Rogers-Ramanujan type, {\em Duke Math. J.} {\bf{32}} (1965), 741-748.

\item[7.] \un{L. J. Slater}, Further identities of Rogers-Ramanujam
type, {\em Proc. London Math. Soc.} (2) {\bf{54}} (1952), 147-167.
\end{enumerate}
\bigskip

Department of Mathematics, 

University of Florida, 

Gainesville, FL 32611, USA
\bigskip

{\it{alladi@math.ufl.edu}}

{\it{alexb@math.ufl.edu}}

\end{document}